# SOME HYNDON'S GENERALIZATIONS TO STARRETT'S METHOD OF SOLVING FIRST ORDER ODES BY LIE GROUP SYMMETRY

Z.M. Mwanzia, K.C. Sogomo

*Zablon M. Mwanzia, Department of Mathematics, P.O Box 536, Egerton.*
*zmusyoka@egerton.ac.ke*

*Dr. K.C. Sogomo, Department of Open and Distance Learning, P.O Box 536, Egerton.*
*kecheso@yahoo.com*

**Abstract**
This paper is centred on solving differential equations by symmetry groups for first order ODEs and is in response to [5]. It also explores the possibility of averting the assumptions in [3] that, in practice finding the solutions of the linearized symmetry condition (the symbols $\xi$ and $\eta$) is usually a much more difficult problem than solving the original ODE but, by inspired guesswork, or geometric intuition, it is possible to ascertain a particular solution of the linearized symmetry condition which will allow the original equation to be integrated. Some differential equations are solved in [5] using Lie group symmetry. Here, the steps used when solving all first order differential equations involve some assumptions and guesses of the form of symmetry for a given differential equation. In [1], there are assertions that many Lie symmetries have tangent vector fields of the form $\xi = \alpha(x)$, $\eta = \beta(x)y + \gamma(x)$, which is a less restrictive ansatz and a more restrictive ansatz (of the form $\xi = c_1 x + c_2 y + c_3$, $\eta = c_4 x + c_5 y + c_6$) can be used to find some common Lie symmetries (for some parameter $\varepsilon$) including translations $(x + \varepsilon, y + \varepsilon)$, scalings $(e^\varepsilon x, e^\varepsilon y)$ and rotations $(x\cos\varepsilon - y\sin\varepsilon, x\sin\varepsilon + y\cos\varepsilon)$. But these ansatz do not always give the required tangent form. Tangent vectors are put in emphasis as there is a one-to-one correspondence between the tangent vectors and Lie group symmetries. The main aim of this paper is to provide a clear computation of tangent vectors corresponding to differential equations by looking for possibility of actually getting the exact form of symmetries as given in [5] using the same (and more) equations. Although Lie group symmetry has been widely employed for solving many equations, there was a neglected problem for first order ODEs. In addition, there is the verification of the connection between prolongation and the linearized symmetry condition as used in solving first order ODEs.

## 1.0 Introduction

The study of symmetry provides one of the most appealing practical applications of group theory. This was extensively shown by Sophus Lie in the nineteenth century [1872–1899]. He investigated the continuous groups of transformations leaving differential equations invariant, creating what is now called the symmetry analysis of differential equations. His aim was to solve non-linear differential equations, which to some extend are may be cumbersome to solve. Intuitively, Lie's method of solving differential equations enables differential equations to be solved in an algebraic approach and as it is put across in [2], symmetry of a differential equation is any transformation of its solution manifold into itself, i.e., a symmetry maps any solution to another solution of the same equation. Furthermore, the symmetry group of a system of differential equations is a group of transformations of independent and dependent variables leaving the set of all solutions invariant. Hence, continuous symmetry transformations are defined topologically and thus are not restricted to local transformations acting on the space of independent and dependent variables and their derivatives. In this sense, any differential equation possesses symmetries.

Lie group symmetry as applied in solving differential equations has been widely employed since its pioneering. This is an overview of first order differential equations which has been influenced by these factors. First, Lie group symmetry as seen in [5] involves assumptions and guesses of the form of symmetry for specific differential equations. Secondly, in [3], it is stated that in practice finding the solutions of the linearized symmetry condition (the symbols $\xi$ and $\eta$) is usually a much more difficult problem than solving the original ODE, but by inspired guesswork, or geometric intuition, it is possible to ascertain a particular solution of the linearized symmetry condition. In general, if the symmetry group of a system of equations is known (or determined), it can be used to generate new solutions from the old ones. Once the tangent vector of a differential equation has been found, the other work of solving it is assumed to easily follow. As indicated, the use of several intuitions and generalizations to find the form of symmetries of first order ODEs as given in [5] will be employed.

## 2.0 Literature Review
### 2.1 Lie groups

As explained in [3], a Lie group is a blending of the algebraic concept of a group and the differential-geometric concept of a manifold. The distinguishing feature of a Lie group from the more general types of groups is that it also carries the structure of a smooth manifold, and the group elements can be continuously varied. Basically, this combination of algebra and calculus leads to a powerful techniques for the study of symmetry.

Definition: Lie group

This is a group $G$ is that is also a smooth manifold and such that the multiplication map
$$\mu : G \times G \to G, \quad \mu : (g,h) \mapsto gh \quad g,h \in G$$
and the inversion map
$$i : G \to G, \quad i : g \mapsto g^{-1} \quad g \in G$$
are smooth maps between manifolds.

An r-parameter Lie group carries the structure of an r-dimensional manifold. A Lie group can also be considered as a topological group (i.e, a group endowed with a topology with respect to which the group operations are continuous) that is also a manifold.

### 2.2 Lie group symmetries

As it is defined in [5], let x $= (x, y)$ and X $= (X, Y)$ be points in the Euclidean plane, and for $\lambda \in \mathbb{R}$, let $P_\lambda : \mathrm{x} \mapsto f(\mathrm{x}, \lambda) = \mathrm{X}$ be a transformation, depending on the parameter $\lambda$, that takes points x to points X. We say the set of transformations $P_\lambda$ is a (additive) transformation group $G$ if the following conditions are satisfied:

1. $P_\lambda$ is one-to-one and onto;
2. $P_{\lambda_2} \circ P_{\lambda_1} = P_{(\lambda_2 + \lambda_1)}$, that is, $f(f(x, \lambda_1), \lambda_2) = f(x, \lambda_2 + \lambda_1)$:
3. $P_0 = I$ (i.e., $f(x, 0) = x$
4. For each $\lambda_1$ there exists a unique $\lambda_2 = -\lambda_1$ such that $P_{\lambda_2} \circ P_{\lambda_1} = P_0 = I$ that is, $f(f(x, \lambda_1), \lambda_2) = f(x, 0) = x$.

If, in addition to these four group properties, $f$ is infinitely differentiable with respect to x and analytic with respect to $\lambda$, we say $G$ is a one-parameter Lie group (or a Lie point transformation). The Lie group symmetries considered are symmetries under a local group, and the group action may not be defined over the whole plane. For instance, the group action

$$P_\lambda: (x, y) \mapsto (X, Y) = \left(\frac{x}{1 - \lambda x}, \frac{y}{1 - \lambda x}\right)$$

is defined only if $\lambda < 1/x$ when $x > 0$ and $\lambda > 1/x$ when $x < 0$.

If $x$ denotes the position of a general point of the object, and if $P: x \mapsto X$ is any symmetry, then we assume that $X$ is infinitely differentiable with respect to $x$. Moreover, since $P^{-1}$ is also a symmetry, $x$ is infinitely differentiable with respect to $X$. Thus $P$ is a ($C^\infty$) diffeomorphism (i.e, a smooth invertible mapping whose inverse is also smooth). Using the above transformation definition and symmetry as defined in [1], a transformation is a symmetry if it satisfies the following:

1. The transformation preserves the structure.
2. The transformation is a diffeomorphism.
3. The transformation maps the object to itself [i.e, a planar object and its image are indistinguishable].

Transformations satisfying 1 and 2 above are symmetries if they also satisfy 3, which is called the symmetry condition. The term "group" is used because the set of symmetries $P_\lambda$ satisfy the axioms of a group, at least for some parameter $\lambda$ sufficiently close to zero.

## 2.3  Lie group symmetries and differential equations
### 2.3.1  Important notations and terms

The notations employed herein have been lifted from [1] and a general first order ODE will be denoted in the following form,

$$\frac{dy}{dx} = h(x, y) \qquad \ldots\ldots\ldots\ldots \quad (1)$$

Invariant solutions of a first order ODE can be characterized by the characteristic equation which is equivalent to

$$Q = \eta - h(x, y)\xi \qquad \ldots\ldots\ldots\ldots \quad (2)$$

All trivial Lie symmetries correspond to $Q = 0$. Actually, the trivial symmetries are useless in solving ODEs. In principle, the nontrivial symmetries can be found by using the method of characteristics and the characteristic equations are

$$\frac{dx}{1} = \frac{dy}{h(x,y)} = \frac{dQ}{h_y(x,y)Q} \qquad \ldots\ldots\ldots\ldots \quad (3)$$

Notably, when seeking for tangent vector fields of a symmetry of a given differential equation, they should satisfy $Q \neq 0$ (should not be trivial) for them to be used to solve the ODE. Hence, any values of $\xi$ and $\eta$ that satisfy $Q \neq 0$ can be used to solve the differential equation.

### 2.3.2  Prolongation and linearized symmetry condition

The solutions of the linearized symmetry condition give the corresponding tangent vectors of a specific differential equation as shown by [5] and [1]. As indicated, in most cases the use of the tangent vectors is common rather than the Lie symmetries themselves, but the tangent vectors are in one-to-one correspondence with the respective Lie symmetries. The same condition can also be expressed in terms of the characteristic equations as seen in [1]. The linearized symmetry condition as used in solving first order ODEs has the following form;

$$\eta_x - \xi_y h^2 + (\eta_y - \xi_x)h - (\xi h_x + \eta h_y) = 0$$

where $\eta, \xi$ and $h$ are as defined.

On the other hand, prolongation of vector fields is another name for the jet of vector fields. An extensive overview and computations for prolongations of vector fields can be found from [3] and [4] among others. For first order ODEs, $h$ as defined, if this equation is invariant under a one-parameter group of transformation, then it can be integrated by quadrature. Taking $G$ to be the group on an open subset $M \subset X \times U \simeq \mathbb{R}^2$ (as defined in [3]), let the infinitesimal generator of the ODE to be of general form,

$$v = \xi(x,y)\frac{\partial}{\partial x} + \eta(x,y)\frac{\partial}{\partial y}$$

From the general prolongation formula, the first prolongation of $v$ is the vector field

$$\mathrm{pr}^{(1)}v = v + \eta^x \frac{\partial}{\partial u_x}$$

where
$$\eta^x = D_x\eta - y_x D_x\xi$$
$$= \eta_x + (\eta_y - \xi_x)y_x - \xi_y y_x^2$$

Utilizing the infinitesimal criterion of invariance, for $v$ to be the infinitesimal generator of the group action, then

$$\mathrm{pr}^{(1)}v(\Delta) = 0 \qquad \qquad \ldots\ldots\ldots\ldots\ldots \quad (4)$$

where the system $\Delta$ as defined here will be $\Delta(x, y, y_x) = y_x - h(x, y)$. $h(x, y)$ is mostly denoted as $h$ for brevity. Assuming that the Jacobian matrix defined here is of rank 1 everywhere and substituting $\Delta$ in condition (4) above,

$$\mathrm{pr}^{(1)}v(\Delta) = \xi\frac{\partial \Delta}{\partial x} + \eta\frac{\partial \Delta}{\partial y} + \eta^x\frac{\partial \Delta}{\partial u_x}$$
$$= -\xi h_x - \eta h_y + \eta_x + (\eta_y - \xi_x)y_x - \xi_y y_x^2 \quad \ldots\ldots \quad (5)$$

From (5) above, it is easily verified that $y_x = h$ and therefore,

$$\mathrm{pr}^{(1)}v(\Delta) = 0$$
$$= -\xi h_x - \eta h_y + \eta_x + (\eta_y - \xi_x)h - \xi_y h^2$$
$$\Rightarrow -(\xi h_x + \eta h_y) + \eta_x + (\eta_y - \xi_x)h - \xi_y h^2 = 0$$
$$\Leftrightarrow \eta_x - \xi_y h^2 + (\eta_y - \xi_x)h - (\xi h_x + \eta h_y) = 0$$

which is the linearized symmetry condition developed.

### 3.0 Symmetries and Standard Methods
As explained in [1], some standard methods of solving ODEs are generalizations of Lie group symmetry methods. Canonical coordinates are associated with a particular Lie group. So all first order ODEs admitting a given one-parameter symmetry group can be reduced to quadrature, using canonical coordinates defined by the group generator.

### 3.1 Symmetries of homogeneous equations
From [3], every ODE of the form

$$\frac{dy}{dx} = F\left(\frac{y}{x}\right), \qquad \qquad \ldots\ldots\ldots\ldots \quad (6)$$

where $F$ depends on the ratio $\left(\frac{y}{x}\right)$,

admits the one-parameter Lie group of scaling symmetries

$$P_\varepsilon: (x, y) \mapsto ((e^\varepsilon x, e^\varepsilon y))$$

The standard solution method for this type of ODE is to introduce new variables

$$r = \frac{y}{x}, \quad s = \ln|x|.$$

These are the canonical coordinates (for $x \neq 0$).

If $F(r) = r$, the symmetries are trivial and the general solution of (6) is $r = c$, that is, $y = cx$. Otherwise, the general solution is,

$$\ln |x| = \int^{\frac{y}{x}} \frac{dr}{F(r) - r} + c$$

In regards to the general linear ODE

$$\frac{dy}{dx} + F(x)y = G(x) \qquad \ldots\ldots\ldots \qquad (7)$$

the homogeneous ODE

$$\frac{du}{dx} + F(x)u = 0$$

is separable and one nonzero solution is

$$u = u_0(x) = \exp\left\{-\int F(x)dx\right\}.$$

We reckon that the principle of linear superposition states that if $y = y(x)$ is a solution of (7), then so is $y = y(x) + \varepsilon u_0(x)$, for each $\varepsilon \in \mathbb{R}$. This principle is equivalent to the statement that (7) has the Lie symmetries

$$P_\varepsilon: (x, y) \mapsto (x, y + \varepsilon u_0(x)).$$

The tangent vector field is,

$$(\xi, \eta) = (0, u_0(x)).$$

## 3.2 Some Hyndon's generalizations of Lie symmetry to Starrett's work

The various steps to be used when solving all first order differential equations in [5] will be presented here. These steps involve assumptions and guesses of the form of symmetry for specific differential equations as illustrated herein. If the symmetry form as found in [5] can be computed for, then the given differential equation can be solved with ease. Here, it is shown and proven that one does not really need to guess the form of symmetry, but the symmetries exactly as given in [5] can actually be found. This will be accomplished by using the more restrictive ansatz and in some cases the less restrictive ansatz and moreover, some situations in which the ansatz may not bear fruits. Moreover, it is evident from [1] that for some first order ODEs, the search for a nontrivial symmetry may be futile although infinitely many such symmetries exist. Here is a review of the steps given in [5], and the latter computations will bring forth the respective symmetry forms.

## 3.3 Steps to solve first-order ODEs

1. Find Lie symmetries of the unknown solutions of the differential equation, say $\frac{dy}{dx} = h(x, y)$, by appealing to the linearized symmetry condition

$$\eta_x - \xi_y h^2 + (\eta_y - \xi_x)h - (\xi h_x + \eta h_y) = 0 \qquad \ldots.. \qquad (8)$$

   To do this, we must make a guess at the form of symmetry.

2. Use the solutions to the linearized symmetry condition (the symbols $\xi$ and $\eta$) to find a coordinate system $(r, s)$ in which the solutions depend on only one of the variables. To do this, we integrate the characteristic equations
   $\frac{dx}{\xi} = \frac{dy}{\eta}$ and $\frac{dx}{\xi} = \frac{dy}{\eta} = ds$ of the orbit.

3. Substitute the canonical coordinates into

$$\frac{ds}{dr} = \frac{s_x + h(x,y)s_y}{r_x + h(x,y)r_y}$$

and solve the differential equation in the canonical coordinate system.

4. Express the solution in the original coordinates.

These steps have no problem except the first step that is of interest here. To avert the problem, the main achievement will be directed towards minimizing the guesses and possibly, find the form of symmetry. To illustrate this, the following examples will be the main feature:

**Example 1.** In solving the differential equation $y' = y^2 x^{-1}$, the symmetry of this equation is not known. But using this ansatz
$$\xi = c_1 y^2 + c_2 x + c_3 y + c_4,$$
$$\eta = c_5 y^2 + c_6 x + c_7 y + c_8.$$

This is a fairly restrictive ansatz. The symmetry condition (8) takes the form,
$$c_6 - (2c_1 y + c_3)y^4 x^{-2} + (2c_5 y + c_7 - c_2)y^2 x^{-1} - \{(c_1 y^2 + c_2 x + c_3 y + c_4)(-y^2 x^{-2})$$
$$+ (c_5 y^2 + c_6 x + c_7 y + c_8)(2yx^{-1})\} = 0$$

Simplifying a little bit and comparing coefficients of powers of $y$, gives:

$y^5: -2c_1 x^{-2} = 0$
$y^4: -c_3 x^{-2} + c_1 x^{-2} = 0$
$y^3: c_3 x^{-2} = 0$
$y^2: c_7 x^{-1} - c_4 x^{-2} - 2c_7 x^{-1} = 0$
$y: -2c_6 - 2c_8 x^{-1} = 0$
$y^0: c_6 = 0$

Simple comparison shows that $c_1 = c_3 = c_6 = c_8 = 0$. Note carefully that $c_2 \neq c_5 \neq 0$ and $c_7 x = c_4$, $\Rightarrow c_7 = c_4 = 0$.

Furthermore, if the expression $c_7 x = c_4$ is to be used, the result will follow that $c_7 = c_4 = 0$.

We are left with only two constants $c_5$ and $c_2$.

The associated vector fields are of the form $\xi = c_2 x$, $\eta = c_5 y^2$ for $c_2 \neq c_5$. It can easily be shown that these are the vector fields of the equation given above.

**Example 2.** Considering the Bernoulli equation $\frac{dy}{dx} = y + y^{-1} e^x$, using the more restrictive ansatz which is
$$\xi = c_1 x + c_2 y + c_3,$$
$$\eta = c_4 x + c_5 y + c_6.$$

The condition (8) with these ansatz will be,
$$c_4 - c_2(y + e^x y^{-1})^2 + (c_5 - c_1)(y + e^x y^{-1}) - \{(c_1 x + c_2 y + c_3)(e^x y^{-1})$$
$$+ (c_4 x + c_5 y + c_6)(1 - e^x y^{-2})\} = 0$$

Following the procedure as in the previous example, we notice that:

$y^2: -c_2 = 0$
$y: -c_1 = 0$
$y^0: c_4 - 3c_2 e^x - c_4 x - c_6 = 0$
$y^{-1}: 2c_5 e^x - c_1 e^x - c_1 x e^x - c_3 e^x = 0$
$y^{-2}: -c_2 e^{2x} + c_4 x e^x + c_6 e^x = 0$

Thus $c_1 = c_2 = 0$.

Using the expressions for $y^{-2}$ and $y$, it follows that $c_1 = c_2 = c_4 = 0$.

The expression for $y^{-1}$ leaves $2c_5 = c_3$ and we are left with only one constant, $c_5$. The associated vector fields are of the form $\xi = c$, $\eta = \frac{1}{2}cy$. Easy computations show that for any constant $c$, these are the vector fields of the Bernoulli equation given above.

### 3.4 Revised steps to solve first order ODEs
The steps to solve a full-fledged first-order ODEs from start to finish now will be as follows:
1. Find Lie symmetries of the unknown solutions of the differential equation $h$, by appealing to the linearized symmetry condition (8). We use a suitable ansatz for the equation to find the respective symmetry.
2. Use the vector fields from (step 1) above (symbols $\xi$ and $\eta$) to find a coordinate system $(r, s)$ in which the solutions depend on only one of the variables. We integrate the characteristic equations
$\frac{dx}{\xi} = \frac{dy}{\eta}$ and $\frac{dx}{\xi} = \frac{dy}{\eta} ds$ of the orbit.
3. Substitute the canonical coordinates into
$$\frac{ds}{dr} = \frac{s_x + h(x,y)s_y}{r_x + h(x,y)r_y}$$
and solve the differential equation in the canonical coordinate system.
4. Express the solution in the original coordinates.

Notably, the solutions of (8) necessarily gives the vector fields (symbols $\xi$ and $\eta$) which are vital to solving the differential equation. As indicated before, once the symbols $\xi$ and $\eta$ have been found, it is easy to solve a differential.

**Example 3:** In solving the Bernoulli equation $\frac{dy}{dx} = y + y^{-1}e^x$ from example 2, the calculated associated vector fields are of the form $\xi = c$, $\eta = \frac{1}{2}cy$.

These latter steps can be seen to follow in sequence as from [5].
The canonical coordinates are found by solving $\frac{dy}{dx} = \frac{\eta}{\xi} = \frac{y}{2}$, to get $r$ and $s$.
Families of functions that remain constant for $r$ are sought here, so $r = c = ye^{-x/2}$.
The second coordinate $s$ is found by integrating $ds = \frac{dx}{c}$ to get $s = x$.
The next step is to express the differential equation in the canonical coordinates by computing
$$\frac{ds}{dr} = \frac{s_x + s_y h}{r_x + r_y h}$$
Proper substitution gives;
$$\frac{ds}{dr} = \frac{1}{-\frac{1}{2}ye^{-x/2} + e^{-\frac{x}{2}}(y + y^{-1}e^x)} = \frac{1}{\frac{1}{2}ye^{-x/2} + y^{-1}e^{x/2}}.$$

Expressing $\frac{1}{2}ye^{-x/2} + y^{-1}e^{x/2}$ in terms of $r$ and $s$ gives the following values
$$\frac{r}{2} + \frac{1}{r},$$
therefore, $\frac{ds}{dr} = \frac{r}{r^2/2 + 1}.$

Which integrates to
$$s = \ln(r^2/2 + 1) + c.$$

Returning to the original coordinates,
$$y = \pm\sqrt{(ce^{2x} - 2e^x)}.$$

## 3.5 Some cases where the ansatz don't work

It may be simple to calculate the symmetries of some differential equation using the more or the less restrictive ansatz. But in some cases, to find the suitable tangent vectors for solving the differential equation using both the ansatz may be fruitless. However, there are some intuitions that can lead to better results and are illustrated using these examples:

**Example 4:** $\frac{dy}{dx} = \frac{y}{x} + x$. This is an easy example it can be easily solved with an integrating factor, but the intention here is to use symmetry methods. Hence, it is vital to start by finding the associated tangent vectors. If the following ansatz is used in finding the Lie symmetries,
$$\xi = c_1 x + c_2 y + c_3$$
$$\eta = c_4 x + c_5 y + c_6$$
It will be noticed that, $c_2 = c_3 = c_6 = 0$ and $2c_1 = c_5$, if the same procedure in deriving the tangent vectors is followed.
This implies that,
$$\xi = c_1 x$$
$$\eta = c_4 x + 2c_1 y$$
These will correspond to the tangent vector as given in [5] if and only if $c_1 = 0$.
In fact, the pair of ansatz
$$\xi = \alpha(x), \quad \eta = \beta(x)y + \gamma(x) \quad \text{and}$$
$$\xi = c_1 x + c_2 y + c_3, \quad \eta = c_4 x + c_5 y + c_6$$
does not seem to help at all in finding the tangent vectors as given by [5], and moreover, it is quite clear that it is difficult to solve the above equation using these derived tangent vector fields. But, the above ODE can be expressed as a homogeneous equation of the form (7) which is separable and one nonzero solution is
$$u = u_0(x) = \exp\left\{-\int F(x)dx\right\}.$$
where $F(x) = -\frac{1}{x}$,
and the solution here is $u = x$.
Hence the tangent vector field is,
$$(\xi, \eta) = (0, x).$$
The method of characteristics can also be used to find the non-trivial symmetries and the characteristic equations are the expression (2).
The characteristic equation here will be $Q = x$.
Since $Q = \eta - h(x,y)\xi$ (where $h$ here is the ODE),
It follows that
$$\eta - x = \xi\left(\frac{y}{x} + x\right) \qquad \ldots\ldots\ldots \quad (9)$$
Clearly, it is easy to compute for the symbols $(\xi, \eta)$ which will satisfy (9).
When $\xi = 0$, $\eta = x$, corresponding to the tangent vector realized earlier.
Similarly when $\xi = -1$, $\eta = -\frac{y}{x}$.
Hence for this ODE, a pair of tangent vectors has been realized,
$$(\xi, \eta) = (0, x)$$

$$(\xi, \eta) = \left(-1, -\frac{y}{x}\right)$$
which can easily be shown to satisfy the symmetry condition for the ODE.
Likewise, $(\xi, \eta) = \left(1, \frac{y}{x}\right)$ also serves the role.

**Example 5:** Look at this differential equation, $\frac{dy}{dx} = \frac{y}{x-y}$.

Following from [5], this equation can be solved by exchanging $x$ and $y$ and computing
$$\frac{dy}{dx} = \frac{D_x x}{D_x y} = \frac{x-y}{y}$$
which is now easily solvable using an integrating factor, but the goal was to solve using Lie group methods.

It does not seem to bring forth the tangent vectors from the ansatz
$$\xi = \alpha(x), \quad \eta = \beta(x)y + \gamma(x) \text{ and}$$
$$\xi = c_1 x + c_2 y + c_3, \quad \eta = c_4 x + c_5 y + c_6.$$
In fact, using the latter one, $c_3 = c_4 = c_6 = 0$ and $c_1 = c_5$ if the procedure in deriving the tangent vectors as used in examples 1 and 2 is followed.
This implies that,
$$\xi = c_1 x + c_2 y$$
$$\eta = c_1 y$$
These will equate to the tangent vector as given in [5] if and only if $c_1 = 0$.
Understandably,
$$\frac{dy}{dx} = \frac{y}{x-y} \quad \Leftrightarrow \quad \frac{dy}{dx} = \frac{\frac{y}{x}}{1-\frac{y}{x}}$$
hence it is wholly dependent on the ratio $\frac{y}{x}$.

From (6), every ODE of the form $\frac{dy}{dx} = F\left(\frac{y}{x}\right)$ where $F$ depends on the ratio $\left(\frac{y}{x}\right)$, admits the one-parameter Lie group of scaling symmetries
$$P_\varepsilon : (x, y) \mapsto ((e^\varepsilon x, e^\varepsilon y))$$
and the corresponding tangent vector is, $(\xi, \eta) = (x, y)$
This can be used to solve the ODE but another symmetry can be found by using the characteristic equation $Q = \eta - h(x,y)\xi$.
Now, $Q = \frac{-y^2}{x-y}$ and
$$y(\xi - y) = \eta(x - y) \qquad \ldots\ldots\ldots \qquad (10)$$
Computing for the symbols $(\xi, \eta)$ which will satisfy (10), easily gives,
when $\xi = y, \eta = 0$,
Similarly when $\xi = x, \eta = y$ corresponding to the tangent vector realized earlier.
Thus, for this ODE we have a pair of tangent vectors
$$(\xi, \eta) = (y, 0)$$
$$(\xi, \eta) = (x, y)$$
which can easily be shown to satisfy the symmetry condition for the ODE and may be used to solve the ODE.

**Example 6:** The following nonlinear ODE does not appear to be solvable by any of the standard methods.

$$\frac{dy}{dx} = \frac{1-y^2}{xy} + 1$$

In finding the point symmetries of this equation, the ansatz as used by [1] is the less restrictive one ($\xi = \alpha(x), \ \eta = \beta(x)y + \gamma(x)$) and here is how the symmetries have been got.
The linearized symmetry condition will be as follows,

$$\beta'y + \gamma' + (\beta - \alpha')\left(\frac{1-y^2}{xy} + 1\right) - \alpha\left(\frac{y^2-1}{x^2 y}\right) + (\beta y + \gamma)\left(\frac{1+y^2}{xy^2}\right) = 0$$

This is a single equation and it can be split into an over determined system of equations by comparing coefficients of powers of $y$, which reveals:

$y^{-2}: \gamma = 0$
$y^{-1}: \frac{\beta - \alpha'}{x} = \frac{\alpha}{x^2} - \frac{\beta}{x}$
$y^{0}: \beta = \alpha'$

The second and third expressions gives, $\alpha' = \frac{\alpha}{x}$.
This ODE is easily solved and the general solution is $\alpha = cx^{-1}$ and therefore $\beta = -cx^{-2}$.
The remaining terms in the linearized symmetry condition provide no further constraints, and hence any tangent vector field of the Lie symmetry group has the form $(\xi, \eta) = (cx^{-1}, -cx^{-2}y)$.

## 4.0 Conclusion & Recommendations

It has been shown that the linearized symmetry condition is actually a case of the prolongation formula truncated for first order ODEs. The same condition can be achieved without using the prolongation formula as seen in [1]. In addition, it has been shown that one can actually find the tangent vectors corresponding to the form of Lie symmetry of specific differential equations. In fact, the tangent vectors of differential equations have been computed for to get the same values exactly as given in [5], and evaded the guess of the form of the symmetries.
The equations that were used herein involve non-linear first order differential equations whose symmetries are achievable. However, there are some cases where the search for the symmetries may not be easy or may fail. In addition, it is also not quite clear how we get symmetries of ODEs which include a one-parameter Lie group of inversions. A perfect example for this is the Riccati equation $\frac{dy}{dx} = \frac{y+1}{x} + \frac{y^2}{x^3}$. The symmetries of this equation include a one-parameter Lie group of inversions. It is not certain about how the symmetries of equations like these are achieved, or maybe they have forms of symmetries which are assumed. Hence, to make the application of Lie group symmetry complete for first order ODEs, there is need to show clearly how such symmetries of this form exist.